\documentclass[11pt]{article}
\usepackage{amsmath,amssymb,amsthm,graphicx,color,bbm}
\usepackage[left=1in,right=1in,top=1in,bottom=1in]{geometry}
\usepackage[british]{babel}
\usepackage{hyperref} 
\usepackage{enumitem}
\usepackage{bm}
\usepackage{csquotes}
\usepackage{tikz}
\usetikzlibrary{decorations.markings}
\usetikzlibrary{automata, positioning,arrows,arrows.meta}
\usetikzlibrary{patterns}
\usetikzlibrary{calc}

\usepackage{pgfplots}
\pgfplotsset{ticks=none}
\pgfplotsset{compat=1.18}

\usepackage[backend=biber,citestyle=numeric,doi=false,isbn=false,url=false,eprint=false,giveninits=true,maxbibnames=10,sortcites]{biblatex}
\AtEveryBibitem{\clearfield{number}}
\AtEveryBibitem{\clearfield{series}}

\DeclareFieldFormat[article,inbook,incollection,inproceedings,patent,thesis,unpublished]{title}{{#1\isdot}}
  \renewbibmacro{in:}{%
  \ifentrytype{article}{}{\printtext{\bibstring{in}\intitlepunct}}}

\makeatletter 
\@mparswitchfalse%
\makeatother
\reversemarginpar
\usepackage[colorinlistoftodos]{todonotes}

\hypersetup{
    colorlinks=false,
    pdfborder={0 0 0},
}

\newtheorem{theorem}{Theorem}[section]

\newtheorem{proposition}[theorem]{Proposition}
\newtheorem{lemma}[theorem]{Lemma}
 \newtheorem{conjecture}[theorem]{Conjecture}

\numberwithin{equation}{section}

\theoremstyle{definition}
\newtheorem{definition}[theorem]{Definition}

\theoremstyle{remark}
\newtheorem{remark}[theorem]{Remark}
\newtheorem{remarks}[theorem]{Remarks}
\newtheorem*{remark*}{Remark}

 \allowdisplaybreaks

\newcommand{\1}[1]{{\mathbbm{1}\mkern -1.5mu}{\{#1\}}}
\newcommand{\2}[1]{{\mathbbm{1}}_{#1}}

\newcommand{\R}{{\mathbb R}}

\newcommand{\Z}{{\mathbb Z}}
\newcommand{\N}{{\mathbb N}}

\newcommand{\ZP}{{\mathbb Z}_+}
\newcommand{\RP}{{\mathbb R}_+}

\DeclareMathOperator{\Exp}{\mathbb{E}}
\let\Pr\relax
\DeclareMathOperator{\Pr}{\mathbb{P}}

\DeclareMathOperator{\Var}{\mathbb{V}ar}

\newcommand{\eps}{\varepsilon}

\newcommand{\re}{{\mathrm{e}}}
\newcommand{\ud}{{\mathrm d}}
\newcommand{\rc}{{\mathrm{c}}}

\newcommand{\cD}{{\mathcal D}}
\newcommand{\cE}{{\mathcal E}}
\newcommand{\cF}{{\mathcal F}}

\newcommand{\cI}{{\mathcal I}}

\newcommand{\DM}{\text{DM}}

\newcommand{\bB}{{\mathbf B}}
\newcommand{\bX}{{\mathbf X}}

\newcommand{\tod}{\overset{\mathrm{d}}{\longrightarrow}}

\newcommand{\Bigmid}{\; \Bigl| \;}

\newcommand{\eqd}{\overset{d}{=}}

\makeatletter
\def\namedlabel#1#2{\begingroup  
    (#2)%
    \def\@currentlabel{#2}%
    \phantomsection\label{#1}\endgroup
}
\makeatother

\newlist{myenumi}{enumerate}{10}
\setlist[myenumi]{leftmargin=0pt, labelindent=\parindent, listparindent=\parindent, labelwidth=0pt, itemindent=!, itemsep=1pt, parsep=4pt}

\newlist{thmenumi}{enumerate}{10}
\setlist[thmenumi]{leftmargin=0pt, labelindent=\parindent, listparindent=\parindent, labelwidth=0pt, itemindent=!}

\title{Energy-constrained Brownian motion}
\author{Oscar Randles\footnote{Department of Mathematical Sciences, Durham University} \and Andrew Wade\footnotemark[1]}
\date{\today}

\begin{document}

\maketitle

\begin{abstract}
For reflected Brownian motion in an interval of length $N$, we consider the first time at which the length of time since the most recent visit to the boundary reaches $M$; this random `lifetime' comes from the interpretation in which the Brownian motion has a finite energy capacity $M$, uses energy at unit rate in the interior, and gets its energy replenished up to $M$ on every visit to the boundary. We prove distributional asymptotics for the lifetime in three different regimes depending on the relative growth of $N$ and $M$; these one-dimensional results parallel results in the random walk case by Wade and Grinfeld. We also present some initial results on the multidimensional  version of the problem, 
for one scaling regime in the case of a ball, 
and conjectures on the behaviour in general domains.
\end{abstract}

\medskip


\medskip


\section{Introduction}
\label{sec:intro}

  The general form of the problem addressed in this note concerns reflected Brownian motion in a smooth, bounded domain in $\R^d$. The quantity of interest is the first time that the current elapsed time since the most recent visit to the boundary exceeds a threshold parameter; this random time we call the \emph{lifetime} of the process (see below for an explanation). 
  A version of this question was studied for random walk in an interval in~\cite{Wade_2023}.
  
  If the domain is scaled by $N$ and the threshold parameter is $M$,
there are three different regimes, as in the random-walk model of~\cite{Wade_2023}, arising from the fact that Brownian motion typically travels distance of order $\sqrt M$ in time $M$. In the \emph{meagre-capacity} regime, $M\ll N^2$, the lifetime ends before an excursion has time to feel the
boundary at a distance. In the \emph{confined-space} regime, $M\gg N^2$, the process typically makes many excursions on the scale of the whole domain before its lifetime is completed on an unusually long excursion. The \emph{critical} regime is $M/N^2\to\rho\in(0,\infty)$, where the energy and spatial scales interact. By Brownian scaling, one can reduce to the single capacity parameter $\rho = M/N^2$ in a fixed domain, with the meagre capacity and confined-space limits being $\rho \to0$ and $\rho \to \infty$, respectively.

In this paper we give a complete picture for the one-dimensional case;
in the meagre capacity regime there is a limit theorem with the so-called  Darling--Mandelbrot law, in the confined-space regime, an explicit normalization gives an exponential limit, and in the critical regime, the limiting law is characterized by its moment generating function. These one-dimensional results (presented in Section~\ref{sec:one-dimensional}) parallel the random-walk results from~\cite{Wade_2023}, but also include a non-asymptotic result for Brownian motion on the half-line, where  a Darling--Mandelbrot statement holds for fixed parameter values. We make modest progress in the higher-dimensional case (Section~\ref{sec:multidim}): for a multidimensional ball, the confined-space limit is proved, with the principal Dirichlet eigenvalue determining the exponential scale.
The method is based on a decomposition of trajectories into high-level and low-level excursions, regeneration, hitting-time estimates, and asymptotic analysis via moment generating functions; the discrete case~\cite{Wade_2023} is similar but does not need the separation into low-level and high-level excursions, which deals with the fact that 
the nontrivial excursions of Brownian motion cannot be listed as a chronological sequence with a first excursion and successive excursions.

We anticipate the results that we obtain here are part of a broader universality picture in which, in the different regimes, the limit law for the lifetime retains more or less data about the underlying domain.   
  At the moment this broader universality picture is speculative (see Section~\ref{sec:universality} below for a discussion). 
  We expect that when $\rho \to 0$ there is a so-called Darling--Mandelbrot limit law with no information about the domain retained, while for $\rho \to \infty$ there is an exponential limit law after a scaling that includes principal Dirichlet eigenvalue and eigenfunction information for the domain. For fixed $\rho \in (0,\infty)$, richer spectral information appears to be relevant.  
  

  We make some comments on the motivation, following~\cite{Wade_2023}. The interpretation as an \emph{energy-constrained Brownian motion} comes from associating  to the diffusing particle a non-negative energy level. Energy decreases at constant rate when the process is in the interior of the domain, and is replenished to maximal level $M$ whenever the process visits the boundary. The first time at which $M$ units of time have elapsed since the last visit to the boundary is thus the first time at which the energy level drops to zero while the process is in the interior, at which point the process terminates, as it no longer has energy to move. This interpretation models a roaming animal that feeds on the boundary of a domain: we refer to~\cite{Wade_2023} for related references in ecology. Similar structure is involved in certain {anticipated rejection} algorithms~\cite{bacher} in computer science, when processes that are taking too long are abandoned.

 \section{One-dimensional  energy-constrained Brownian motion}
 \label{sec:one-dimensional}

\subsection{Overview}

Write $I_\infty := \RP :=[0,\infty)$ and set $I_N := [0,N]$ for $N \in (0, \infty)$. 
Also set $\partial I_\infty := \{ 0\}$ and $\partial I_N := \{ 0, N\}$ for $N < \infty$.
For $N \in  (0,\infty) \cup \{ \infty\}$ we consider reflected Brownian motion on $I_N$,
reflected at $\partial I_N$, which we denote by $X = (X_t )_{t \in \RP}$. 
We suppose that $X_0 = 0$.

For $t \in \RP$, define the time elapsed since the most recent boundary visit
\[ e_t := t - \sup \{ s \in [0,t] : X_s \in \partial I_N \} .\]
Then for $M \in (0,\infty)$, define
\begin{equation}
    \label{eq:lambda-def}
    \lambda := \inf \{ t \in \RP : e_t \geq M \} .
\end{equation}
The interpretation of $\lambda$ is as the lifetime of an energy-constrained Brownian motion,
where the energy of the Brownian motion starts at its maximum level $M$, energy decreases at unit rate while $X_t \in I_N \setminus \partial I_N$, and whenever $X_t \in \partial I_N$ energy is reset to level $M$: $\lambda$ is then the first time that the energy falls to level~$0$. With this picture in mind, we say that $X$ \emph{terminates} at time $\lambda$.

The goal of this section is to present distributional convergence results for $\lambda$, as $M , N \to \infty$. We expect that an alternative proof strategy, using It\^o excursion theory, might lead to shorter proofs, but our approach, decomposing into low-level and high-level excursions, is relatively elementary.
The first case that we consider,
in Section~\ref{sec:BM-on-RP}, is  the meagre capacity regime.

\subsection{The meagre-capacity regime}
\label{sec:BM-on-RP}

To state our result in this case, we need the following definition.

\begin{definition}\label{def: Definition 1.3}
Set
\[ \cI (t) := t \int_0^1 u^{-1/2} \re^{ut} \ud u .\]
 We say $\xi$ has the Darling--Mandelbrot~\cite{darling,lew,lh,bacher,louchard} distribution with parameter $1/2$,
   and write  $\xi\sim \DM(1/2)$, if 
    \begin{equation*}
        \Exp \bigl[ \re^{t \xi} \bigr]=\frac{1}{\re^t-\cI(t)}\qquad \forall t\in \R \text{ with }t<t_0,
    \end{equation*}
    where $t_0:=\inf\{t\in \R:\re^t-\cI(t)\leq 0\} \approx 0.8540326566$ (see Lemma A.1 in~\cite{Wade_2023}).
    \end{definition}

\begin{theorem}
    \label{thm:meagre}
    Let $M, N \in (0,\infty)$ with $M/ N^2  \to 0$. Let $X$ be reflected Brownian motion on $I_N$ with lifetime $\lambda$ as defined at~\eqref{eq:lambda-def}.  Then,
        for $\xi\sim \DM(1/2)$, as $M \to \infty$, 
        \[ \frac{\lambda}{M} \tod 1 +\xi.\]
\end{theorem}
 \begin{remark}
        The $\DM(1/2)$ limit is exactly the same as appears in the corresponding random-walk result, Theorem 2.3 of~\cite{Wade_2023}. In fact, we conjecture that the $\DM(1/2)$ limit appears also in the multidimensional meagre-capacity setting, as we explain in Section~\ref{sec:universality} below.
    \end{remark}

In the meagre-capacity case, the intuition is that Brownian motion
finds it hard to traverse distance $N$ in time $M$, and so 
we do not feel the presence of the other boundary. Thus a good approximation is to consider the case where  $X$ is reflected Brownian motion on $\RP$, i.e., we can write $X_t = | B_t |$
where $B := (B_t)_{t \in \RP}$ is standard Brownian motion on $\R$.

    Here is our main result for that case, from which Theorem~\ref{thm:meagre} will follow. Note that (due to Brownian scaling, as we will show) we obtain a distributional equality rather than just a limit statement in this case.

    \begin{theorem}\label{thm:half-line}
        Let $M \in (0,\infty)$. 
    Let $X$ be reflected Brownian motion on $I_\infty$ with lifetime $\lambda$ as defined at~\eqref{eq:lambda-def}. 
 Then,
        for $\xi\sim \DM(1/2)$,
        \begin{equation*}
             \frac{\lambda}{M} \eqd 1+\xi.
        \end{equation*}
    \end{theorem}
   
    In outline, our proof of Theorem~\ref{thm:half-line} has large parts in common
    with the proof of the corresponding random-walk result, Theorem 2.3 from~\cite{Wade_2023}. The latter
    proof used an excursion decomposition of the random walk. In the case of Brownian motion,
    the excursions cannot be simply ordered in the same way, so a different idea is needed.
    We proceed by splitting the trajectory in a different
    way, as we now describe.

Throughout we write $\N :=\{1,2,3,\ldots\}$ and $\ZP:= \{0\} \cup \N$. 
Recall that $X_0 = 0$. Define
$\gamma_0 := 0$ and then, recursively, for $n \in \N$,
\begin{align}
    \label{eq:sigma-def}
\sigma_n & := \inf \{ t \geq \gamma_{n-1} : X_t = 1 \} ; \\
    \label{eq:gamma-def}
    \gamma_n & := \inf \{ t \geq \sigma_n : X_t  = 0 \}.
\end{align}
Then $0=\gamma_0 < \sigma_1 < \gamma_1 < \sigma_2 < \cdots$ are alternating
sequences of stopping times. 
We say that over time interval $[ \sigma_n, \gamma_n]$ the process is
undertaking a \emph{high-level excursion} while over
$[\gamma_n, \sigma_{n+1}]$ it is undertaking a \emph{low-level} excursion;
the terminology reflects the fact that $X_t \in [0,1]$ for all $ t \in [\gamma_n, \sigma_{n+1}]$.

Define
\begin{equation}
    \label{eq:kapp-def}
    \kappa :=\max \{n\in \ZP: \gamma_n \leq \lambda\},
\end{equation}
the number of completed high-level  excursions
before termination. 
Denote the duration of the $n$th low-level excursion by 
    \begin{equation}\label{eq:2.7}
         \ell_n :=\sigma_n -\gamma_{n-1}  ,
    \end{equation}
and the duration of the $n$th high-level excursion by
    \begin{equation}\label{eq:2.8}
        \nu_n :=\gamma_n -\sigma_n .
    \end{equation}

We show that the durations of low-level excursions are tightly controlled; note that, by the Markov property, every $\ell_n$ defined by~\eqref{eq:2.7} has the same distribution as $\sigma_1$, the first hitting time of~$1$ by reflected Brownian motion started at~$0$, equivalently,
$\sigma_1 = \inf \{ t \in \RP : |B_t | =1\}$.

\begin{lemma}
\label{lem: Lemma 2.1}
For constants $C < \infty$ and $c >0$ it holds that
 $\Pr (\sigma_1>t)\leq C \re^{-c t }$ for all $t \in \RP$.
\end{lemma}
\begin{proof}
Let $\cF_t = \sigma (B_s : 0 \leq s \leq t)$. Note that, for $n \in \ZP$,
\begin{align*}
\Pr ( \sigma_1 > n+1 \mid \cF_n) & = \2{ \{\sigma_1 > n\} } \Pr ( \sigma_1 > n+1 \mid \cF_n ) \\
& \leq \2{ \{\sigma_1 > n\} } \Pr ( | B_{n+1} | < 1 \mid \cF_n ) .
\end{align*}
On the event $\{ \sigma_1 > n \}$ we have $B_n \in (-1,1)$. If $B_n \in (-1,0]$ then for $Z = B_{n+1} - B_{n}$ standard Gaussian independent of $\cF_n$, 
$\Pr ( | B_{n+1} | < 1 \mid \cF_n ) \leq \Pr ( Z > -1 \mid \cF_n ) = \Phi (1) =: p$,
where $\Phi$ is the standard Gaussian distribution function. Similarly, if $B_n \in [0,1)$ then for $Z = B_{n+1} - B_{n}$ standard Gaussian, 
$\Pr ( | B_{n+1} | < 1 \mid \cF_n ) \leq \Pr ( Z < 1 \mid \cF_n ) = p$. So we conclude that
\[ \Pr ( \sigma_1 > n+1  ) \leq p \cdot \Pr ( \sigma_1 > n ) .\]
Hence $\Pr (\sigma_1 > n) \leq p^n \leq \re^{- c n}$, where $c := \log (1/p) >0$.
For general $t \geq 0$, $\Pr (\sigma_1 > t ) \leq \Pr (\sigma_1 > \lfloor t \rfloor )$, and the result follows with $C := \re^c$.
\end{proof}
\noindent
 
Since low-level excursions are brief, we will show that   $\lambda$, at the scale of interest, is well-approximated by quantities defined purely in terms of high-level excursions.
Observe that on a high-level excursion we do not see the reflections
of the reflected Brownian motion, so each $\nu_n$ as defined at~\eqref{eq:2.8}
has the same distribution as the hitting time of~$0$ by standard Brownian motion on $\R$ started at $B_0 = 1$. That is, $\nu_n$ has a positive ($1/2$)-stable (L\'evy) distribution: from the
reflection principle for Brownian motion~\cite[pp.~173--5]{fellerv2}
\begin{equation}
\label{eq:nu-def}
\Pr ( \nu_1 > t ) = \Pr \bigl( \sup_{0 \leq s \leq t} B_s < 1 \bigr) = 2 \Phi ( t^{-1/2} ) - 1 
.
\end{equation} 
Define 
\begin{equation}
    \label{eq:kappa-high}
\kappa_1 := \min \{ n \in \N : \nu_n > M \} -1, \text{ and }
\theta (M) := \Pr ( \nu_1 > M),
\end{equation}
so that $\kappa_1$ is the number of high level excursions before the first one exceeding duration~$M$, an event of probability $\theta (M)$.
Observe that $\nu_n >M$ implies $\lambda < \gamma_n$ and hence $\kappa < n$; hence
\begin{equation}
    \label{eq:kappa-kappa-high}
    \kappa \leq \kappa_1, \text{ a.s.};
\end{equation}
there might be fewer completed high-level excursions if termination occurs on a low-level excursion, for example. Fix $\eps>0$ and define
\begin{equation}
    \label{eq:kappa-high-new}
    \kappa_2 := \min \{ n \in \N : \nu_n > M - M^\eps \} - 1.
\end{equation}
Comparison with~\eqref{eq:kappa-high} shows that $\kappa_2 \leq \kappa_1$, a.s.
We will approximate $\lambda$ by the two quantities
\begin{equation}
    \label{eq:lambda-high}
    \Lambda_1 := M+\sum_{n=1}^{\kappa_1} \nu_n, \text{ and } \Lambda_2 := M + \sum_{n=1}^{\kappa_2} \nu_n.
\end{equation}
We also define
\[ L_M := \sum_{n=1}^{\kappa_1+1} \ell_n .\]
Define the event 
   \begin{equation}\label{eq:2.11}
        E_M :=\{\kappa_1 <M^{1/2+\varepsilon}, \, \ell_n<M^\varepsilon \;\forall n\leq \kappa_1+1 \}.
    \end{equation}
    The next result quantifies how $\Lambda_1$ and $\Lambda_2$ serve
    as approximations to $\lambda$.

    \begin{lemma}
    \label{lem:lambda-approx} Let $\eps \in (0,1)$. 
Define $\Lambda_1$ and $\Lambda_2$ as at~\eqref{eq:lambda-high}, and consider event $E_M$ defined at~\eqref{eq:2.11}. Then (i) $\lambda \leq \Lambda_1 + L_M$, and (ii) on $E_M$, $\lambda \geq \Lambda_2$.
    \end{lemma}
    \begin{proof}
    Note that $n = \kappa_1+1$ is the first high-level excursion with $\nu_n > M$. Whatever has happened
    earlier, termination must occur no later than $M$ units after $\sigma_{\kappa_1+1}$, and hence
    \[ \lambda \leq \sigma_{\kappa_1+1} + M = M + \sum_{j=1}^{\kappa_1+1} \ell_j + \sum_{j=1}^{\kappa_1} \nu_j = \Lambda_1 + L_M ,\]
proving (i). On $E_M$, every high-level excursion is preceded by a low-level excursion of duration less than $M^\eps$, meaning that $\gamma_n - \gamma_{n-1} < M$ for all $n \leq \kappa_2$. Thus termination cannot occur before time $\gamma_{\kappa_2} + M$, and so, in particular,
        \[ \lambda \geq \gamma_{\kappa_2} + M \geq M + \sum_{n=1}^{\kappa_2} \nu_n = \Lambda_2 .\]
        This completes the proof.
    \end{proof}

To study $\Lambda_1$ and $\Lambda_2$, we define the moment generating function of the conditional distribution
\begin{equation}
\label{eq:psi-def}
    \psi_M(s):=\Exp\left[\re^{s\nu_1}\mid \nu_1\leq M\right].
\end{equation}
Recalling Definition~\ref{def: Definition 1.3}, set $K(t):= \re^t-\cI(t)$.

\begin{lemma}\label{lem: Lemma 2.2}
    Consider an energy-constrained Brownian motion on $\RP$. Then, as $M \to \infty$,
    \begin{equation}
    \label{eq:theta-asymp}
        \theta(M)\sim \sqrt{\frac{2}{\pi}}M^{-1/2},
    \end{equation}
    and, uniformly for $t$ in any  compact $R \subset \R$, 
    \begin{equation}
    \label{eq:psi-asymp}
        \lim_{M\to \infty}\left[\left(\psi_M(t/M)-1\right)M^{1/2}\right]=\sqrt{\frac{2}{\pi}}\left(1-K(t)\right).
    \end{equation}
\end{lemma}
\begin{proof}
The proof is similar to that of Proposition~4.2 in~\cite{Wade_2023}, using the Brownian motion estimates in place of similar estimates for random-walk return times. First observe that
\[
 \sqrt{\frac{2}{\pi}} t^{-1/2} -   \sqrt{\frac{1}{2\pi}} t^{-3/2}  \leq
\sqrt{\frac{2}{\pi}} t^{-1/2} \re^{-1/(2t)} \leq
2 \Phi ( t^{-1/2} ) -1 = \sqrt{\frac{2}{\pi}} \int_0^{t^{-1/2}} \re^{-z^2/2} \ud z 
\leq \sqrt{\frac{2}{\pi}} t^{-1/2}
,\]
which with~\eqref{eq:nu-def} gives~\eqref{eq:theta-asymp}.
For the second claim we use a suitable integration by parts formula (see Lemma 4.3 from~\cite{Wade_2023}) to get
\begin{equation*}
    \Exp\left[(e^{\nu_kt/M}-1)\2{\{\nu_k\leq M\}}\right]=(1-\re^t)\Pr(\nu_k>M)+t\int_0^1\re^{tu}\Pr(\nu_k>uM) \ud u.
\end{equation*}
For convenience set $c: = \sqrt{2/\pi}$ and $q_M (u) := M^{1/2} \Pr ( \nu_1 > u M)$, $u >0$.
Then, as already proved, $0 \leq \Pr (\nu_1 > s) \leq c s^{-1/2}$ and hence $0 \leq q_M (u) \leq c u^{-1/2}$, while for each fixed $u >0$, $q_M(u) \to c u^{-1/2}$ as $M \to \infty$. Since $u^{-1/2}$ is integrable on $(0,1)$, dominated convergence gives that
\[ \lim_{M \to \infty} \int_0^1 | q_M (u) - c u^{-1/2} | \ud u = 0 .\]
If $R \subset \R$ is compact and $T = \sup_{t \in R} | t|$, then
\[ \sup_{t \in R} \left| M^{1/2} t \int_0^1 \re^{tu} \Pr ( \nu_1 > u M ) \ud u - c t \int_0^1 \re^{tu} u^{-1/2} \ud u \right| \leq T \re^T \int_0^1 | q_M (u) - c u^{-1/2} | \ud u  \to 0.\]
Thus uniformly for $t \in \R$,
\[ M^{1/2} t \int_0^1 \re^{tu} \Pr ( \nu_1 > u M ) \ud u \to c \cI(t).\]
Hence, using~\eqref{eq:theta-asymp},
\[ M^{1/2} \Exp\left[(\re^{\nu_kt/M}-1)\2{\{\nu_k\leq M\}}\right]
\to 
c \bigl( 1 -\re^t + \cI (t) \bigr) = c ( 1 - K(t)) .\]
Then
\begin{equation*}
    \psi_M(t/M) - 1 = \Exp\left[ \re^{\nu_kt/M}-1 \Bigmid \nu_k\leq  M\right]= \frac{\Exp\left[(\re^{\nu_kt/M}-1)\2{\{\nu_k\leq M\}}\right]}{\Pr(\nu_k\leq M)}
\end{equation*}
and $\Pr(\nu_k\leq M)\to 1$ as $M\to\infty$ (not depending on~$t$), from which we obtain~\eqref{eq:psi-asymp}.
\end{proof}

Next we show that the desired asymptotics for $\lambda$ are achieved by
both $\Lambda_1$ and $\Lambda_2$.

\begin{lemma}\label{lem: Lemma 2.3}
Recall the definitions  of $\Lambda_1$ and $\Lambda_2$ from~\eqref{eq:lambda-high}. 
   With $\xi\sim \DM(1/2)$, 
    as $M \to \infty$, 
    \begin{equation*}
        \frac{\Lambda_1}{M} \tod 1+\xi, \text{ and } \frac{\Lambda_2}{M} \tod 1+\xi.
    \end{equation*}
\end{lemma}
\begin{proof}
Given the estimates in Lemma~\ref{lem: Lemma 2.2}, the proof is similar to
the proof of Theorem 2.3 in~\cite{Wade_2023}. 
Regeneration and the definitions at~\eqref{eq:kappa-high} imply that
 $\Pr(\kappa_1 =k)=\theta(M)(1-\theta(M))^k$ for $k\in \ZP$. Then 
    \begin{equation*}
        \Exp[r^{\kappa_1}]=\sum_{k=0}^\infty r^k(1-\theta(M))^k\theta(M)=\frac{\theta(M)}{1-(1-\theta(M))r}, \text{ for }|(1-\theta(M))r|<1.
    \end{equation*}
    Conditional on $\kappa_1$ the first $\kappa_1$ high-level excursions  are independent and each have the conditional law described by~\eqref{eq:psi-def}. 
    Write $H_{\kappa_1} = \sum_{k=1}^{\kappa_1} \nu_k$. 
    By conditional independence,
    \begin{align*}
        \Exp[\re^{sH_{\kappa_1}}] & = \Exp[(\psi_M(s))^{\kappa_1}]    =\frac{\theta(M)}{1-(1-\theta(M))\psi_M(s)},  \text{ for } |(1-\theta(M))\psi_M(s)|<1.
    \end{align*}
   By Lemma~\ref{lem: Lemma 2.2},  
   \begin{equation*}
       (1-\theta(M))\psi_M(t/M)=1-M^{-1/2}\sqrt{\frac{2}{\pi}}K(t)+o(M^{-1/2})
   \end{equation*}
   meaning $|(1-\theta(M))\psi_M(t/M)|<1$ holds for each fixed $t < t_0$ and $M$ sufficiently large. Hence,
   \begin{equation*}
       \lim_{M\to\infty}\Exp[\re^{tH_{\kappa_1}/M}]=\lim_{M\to\infty}\sqrt{\frac{\pi}{2}}\frac{\theta(M)}{M^{-1/2}K(t)+o(M^{-1/2})}=\frac{1}{K(t)}, \text{ for all } t < t_0.
   \end{equation*}
   Since, by~\eqref{eq:lambda-high}, $\Lambda_1 = M+H_{\kappa_1}$, and recalling
   Definition~\ref{def: Definition 1.3}, the result for $\Lambda_1$ follows.
   The case of $\Lambda_2$ follows the same lines, with essentially the same asymptotic estimates
   appearing when one replaces $\kappa_1$ by $\kappa_2$.
\end{proof}

\begin{proof}[Proof of Theorem~\ref{thm:half-line}]
Fix $\eps \in (0,1/4)$.
First we show that, with $E_M$ as defined at~\eqref{eq:2.11}, $\Pr ( E_M) \to 1$.
Note that
    \begin{equation*}
        \Pr(E_M^\rc)\leq \Pr(\kappa_1 \geq M^{1/2+\varepsilon})+\Pr\left(\kappa_1 \leq M^{1/2+\varepsilon},\bigcup_{k=1}^{\kappa_1+1} \{\ell_k>M^\varepsilon\}\right).
    \end{equation*}
    Here 
    \begin{equation*}
        \Pr(\kappa_1 \geq M^{1/2+\varepsilon})\leq (1-\theta(M))^{M^{1/2 +\varepsilon}}=\exp{\{M^{1/2+\varepsilon}\log(1-\theta(M))\}} \leq  \exp\{-M^{\varepsilon/2} \} ,
    \end{equation*}
    for large $M$, by~\eqref{eq:theta-asymp}. So $\Pr(\kappa_1 \geq M^{1/2+\varepsilon})\to 0 $ as $M\to\infty$. Moreover, Lemma~\ref{lem: Lemma 2.1} shows that
$\Pr(\ell_n>M^\varepsilon)\leq C \exp(-c M^{\varepsilon})$, and so a union bound shows
$\Pr ( E_M ^\rc ) \to 0$.

Now by statement (i) in Lemma~\ref{lem:lambda-approx}, for every $y \in \RP$,
    \begin{equation*}
        \Pr\left(\frac{\lambda}{M}\leq y\right)\geq \Pr\left(\frac{\Lambda_1+L_M}{M}\leq y,E_M \right)\geq   \Pr\left(\frac{\Lambda_1}{M}+\frac{(M^{1/2+\varepsilon}+1)M^\eps}{M}\leq y\right)-\Pr\left(E_M ^\rc\right), 
    \end{equation*}
    which converges to $\Pr ( 1 + \xi \leq y)$, by the $\Lambda_1$ part of Lemma~\ref{lem: Lemma 2.3}
    and the fact that $\Pr ( E_M ^\rc ) \to 0$. On the other hand,
    by statement (ii) in Lemma~\ref{lem:lambda-approx}, for every $y \in \RP$,
    \begin{equation*}
    \Pr\left(\frac{\lambda}{M}\leq y\right)
    \leq \Pr\left(\frac{\Lambda_2}{M}\leq y,E_M \right)+\Pr\left(E_M^\rc\right) ,
    \end{equation*}
     which also converges to $\Pr ( 1 + \xi \leq y)$, by the $\Lambda_2$ part of Lemma~\ref{lem: Lemma 2.3}.    Thus we conclude that $\lambda /M$ converges in distribution to $1+\xi$.

     Finally, to emphasize the dependence on $M$ write $\lambda^{(M)}$ instead of $\lambda$. By Brownian scaling,
       $M^{-1/2}X_{Mt}$ is again reflected Brownian motion on $\RP$, and its elapsed-since-boundary time is $M^{-1}\ell_{Mt}$, hence $\lambda^{(M)}/M$ has the same distribution as $\lambda^{(1)}$.
       Since the preceding argument shows that this distribution converges, as $M\to\infty$, to that of $1+\xi$, it must already equal the distribution of $1+\xi$.
\end{proof}

Finally, we can deduce the meagre-capacity result, Theorem~\ref{thm:meagre}.

\begin{proof}[Proof of Theorem~\ref{thm:meagre}.]
Couple reflected Brownian motion on $[0,N]$, started from $0$, with reflected
Brownian motion on $\RP$, started from $0$; the paths of the two processes
agree until the first time one (hence both) processes hit $N$. Thus, writing $\lambda$ for the lifetime of the process on $[0,N]$,
and using that the lifetime of the process on $\RP$ satisfies the distributional
identity in Theorem~\ref{thm:half-line},
\[ \left| \Pr \left( \frac{\lambda}{M} \leq t \right) - \Pr ( 1 + \xi \leq t ) \right| \leq \Pr \left( \sup_{0 \leq s \leq t M} | B_s | \geq N \right) .
 \]
 Here, by Brownian scaling and a union bound using~\eqref{eq:nu-def},
 \[
 \Pr \left( \sup_{0 \leq s \leq t M} | B_s | \geq N \right)
 = \Pr \left( \sup_{0 \leq s \leq t M/N^2} | B_s | \geq 1 \right)
 \leq 4 \left[ 1 - \Phi ( N t^{-1/2} M^{-1/2} )\right] \to 0 ,
 \]
 for each fixed $t \in \RP$, since $ N M^{-1/2} \to \infty$, by hypothesis.
\end{proof}

\subsection{The confined-space regime}
\label{sec:confined-space}

In the confined space limit, the intuition is that the process will traverse the whole interval many times, and termination occurs on the rare event of an unusually long excursion. In such meta-stable situations it is common to observe an exponential distribution in the limit after suitable scaling~\cite{aldous}, and the next result is of that type.

    \begin{theorem}\label{thm:confined-space}
         Let $M, N \in (0,\infty)$ with $M/ N^2  \to \infty$. Let $X$ be reflected Brownian motion on $I_N$ with lifetime $\lambda$ as defined at~\eqref{eq:lambda-def}.  Then, with $\cE$ a unit-mean exponential random variable, as $M \to \infty$, 
        \begin{equation*}
            \frac{4\lambda}{N^2}\exp\left\{ - \frac{\pi^2 M}{2N^2}\right\} \tod \cE.
        \end{equation*}
    \end{theorem}

Let us make the dependence on $N,M$ explicit and write $\lambda_{N,M} = \lambda$
for the $\lambda$ with spatial domain $I_N$ and energy capacity $M$
defined through~\eqref{eq:lambda-def}; Theorem~\ref{thm:confined-space} 
concerns asymptotics for $\lambda_{N,M}$ as $N, M \to \infty$. Brownian scaling shows that $\lambda_{N,M}$ has the same distribution as $N^2 \lambda_{1,\rho}$ where $\rho := M/N^2$. Thus in the remainder of this section we will work with $\lambda_{1,\rho}$ where $\rho \to \infty$.

We modify the high-level and low-level excursions construction from the previous section, because now we expect many more excursions and we still want to make sure the time spent in low-level excursions is negligible. Recall that $X_0 = 0$. Define
$\gamma_0 := 0$ and then, recursively, for $n \in \N$,
\begin{align}
    \label{eq:sigma-def2}
\sigma_n & := \inf \{ t \geq \gamma_{n-1} : X_t  \in \{ \rho^{-1} , 1 - \rho^{-1} \}\} ; \\
    \label{eq:gamma-def2}
    \gamma_n & := \inf \{ t \geq \sigma_n : X_t  \in \{0,1\} \}.
\end{align}
Again, we say that over time interval $[ \sigma_n, \gamma_n]$ the process is
undertaking a {high-level excursion} while over
$[\gamma_n, \sigma_{n+1}]$ it is undertaking a {low-level} excursion; during low level excursions, the process is within distance $\rho^{-1}$ of the boundary.
We use the same notation as at~\eqref{eq:kapp-def}--\eqref{eq:2.8}, namely
\[    \kappa :=\max \{n\in \ZP: \gamma_n \leq \lambda\}, 
\]
the number of completed high-level   excursions
before termination, 
and
  \[ 
         \ell_n :=\sigma_n -\gamma_{n-1}  , \text{ and } 
        \nu_n :=\gamma_n -\sigma_n ,
\]
the durations of 
low-level and high-level excursions, respectively.
Similarly to~\eqref{eq:kappa-high}, define
\[ \kappa_1 := \min \{ n \in \N : \nu_n > \rho \} -1, \text{ and }
\theta (\rho) := \Pr ( \nu_1 > \rho), \]
so that $\kappa_1$ is the number of high level excursions before the first one exceeding duration~$\rho$, an event of probability $\theta (\rho)$.
Also, similarly to~\eqref{eq:kappa-high-new}, define
\[
\kappa_2 := \min \{ n \in \N : \nu_n > \rho - \rho^{-1/2} \} - 1,
\]
and then, similarly to~\eqref{eq:lambda-high}, define
\[
    \Lambda_1 := \rho +\sum_{n=1}^{\kappa_1} \nu_n, \text{ and } \Lambda_2 := \rho + \sum_{n=1}^{\kappa_2} \nu_n.
\]
We will use $\Lambda_1$ and $\Lambda_2$ to obtain bounds on $\lambda$ in a similar way to the previous section.

\begin{lemma}
\label{lem:lambda-asymptotics-confined}
With the above definitions, for $\cE$ a unit-mean exponential, it holds that, as $\rho \to \infty$,
\[ 4 \re^{- {\pi^2} \rho / {2} } \Lambda_1 \tod \cE, \text{ and }
4 \re^{- {\pi^2} \rho / {2} } \Lambda_2 \tod \cE .\]
\end{lemma}

This result will be proved by an application of
the following exponential convergence result, which 
is Lemma 4.6
in~\cite{Wade_2023}
and can
be deduced from e.g. Theorem 3.2.4 of~\cite[p.~85]{kalashnikov}.

\begin{lemma}\label{lem:exponential}
Let $K_\rho\in \Z_+$ satisfy $\Pr(K_\rho=k)=(1-p_\rho)^kp_\rho$ for $k\in \Z_+$, where $p_\rho\in (0,1)$ and $\lim_{\rho \to\infty}p_\rho=0$. Suppose also that $Y_\rho,Y_{\rho,1},Y_{\rho,2}\ldots $ are i.i.d., $\RP$-valued and independent of $K_\rho$, with $\Exp[Y^2_\rho]=\sigma_\rho^2$ and $\Exp[Y_\rho]=\mu_\rho>0$. Let $Z_\rho:=\sum_{i=1}^{K_\rho}Y_{\rho,i}$. Assuming that 
\begin{equation*}
    \lim_{\rho\to\infty }\frac{\sigma_\rho^2p_\rho}{\mu_\rho^2}=0,
\end{equation*}
it is the case that, as $\rho \to\infty$,
\begin{equation*}
    \frac{p_\rho Z_\rho}{\mu_\rho} \tod \mathcal{E},
\end{equation*}
where $\mathcal{E}$ is exponentially distributed with parameter~$1$.    
\end{lemma}

To apply Lemma~\ref{lem:exponential} to obtain Lemma~\ref{lem:lambda-asymptotics-confined} we will need  
\begin{equation}
\label{eq:rho-full}
\mathbb P(\nu_1>s)
=\frac4\pi\sum_{j=0}^{\infty}
\frac{\sin((2j+1)\pi/\rho)}{2j+1}
\re^{-(2j+1)^2 \pi^2 s/2},
\end{equation}
which follows from formula (7.15) of~\cite{morters2010brownian}. 
In particular, we have the asymptotics
\begin{equation}
\label{eq:rho-tail}
    \Pr ( \nu_1 > \rho ) = \frac{4}{\rho}(1+o(1))\exp\left\{-\frac{\pi^2 \rho}{2}\right\},
\end{equation}
which can be found in~\cite[p~342]{feller1968}. We will also need that 
\begin{equation}
    \label{eq:rho-moments}
    \Exp ( \nu_1 ) \sim \frac{1}{\rho}, \text{ and } \Var ( \nu_1 ) \sim \frac{1}{3\rho} ;
\end{equation}
the statement for $\Exp (\nu_1)$ can be found in
Theorem 2.49 of~\cite{morters2010brownian}, the variance is Exercise 2.15(c) in the same reference.

\begin{proof}[Proof of Lemma~\ref{lem:lambda-asymptotics-confined}.]
    First observe that, for $\beta >0$,
    \begin{align*}
        \Exp ( \nu_1^\beta \1 { \nu_1 \geq \rho } ) & = \int_0^\infty \Pr ( \nu_1^\beta \1{ \nu_1 > \rho} > y ) \ud y \\
        & \leq \rho^\beta \Pr ( \nu_1 \geq \rho ) + \beta \int_\rho^\infty s^{\beta-1} \Pr ( \nu_1 > s ) \ud s.
    \end{align*}
    Next note that from~\eqref{eq:rho-full},
using $|\sin u|\le |u|$, we obtain, uniformly for $\rho>2$ and $s\ge1$,
\[
\mathbb P(\nu_1>s)\le \frac C\rho \re^{-\pi^2 s/2}.
\]
    Then we get, for all $\rho$ large enough, for $C < \infty$,
    \begin{align*}
        \Exp ( \nu_1^\beta \1 { \nu_1 \geq \rho } ) &
        \leq C \rho^{\beta-1} \re^{-\pi^2 \rho /2}
        + C \rho^{-1} \int_\rho^\infty s^{\beta-1} \re^{-\pi^2 s /2} \ud s\\
        & \leq C \rho^{\beta-1} \re^{-\pi^2 \rho /2}
        + C \rho^{-1} \re^{-\pi^2 \rho /4} \int_1^\infty s^{\beta-1} \re^{-\pi^2 s /4} \ud s \\
        & \leq C  \rho^{\beta-1} \re^{-\pi^2 \rho /4} .
        \end{align*}
        Then we have
        \begin{align*}
        \Exp ( \nu_1 \mid \nu_1 < \rho ) \Pr ( \nu_1 < \rho )
        & = \Exp (\nu_1 ) - \Exp ( \nu_1 \1 { \nu_1 \geq \rho } ) ,\end{align*}
        which together with~\eqref{eq:rho-moments} and the fact that
        $\Pr ( \nu_1 < \rho ) \to 1$ yields $\Exp ( \nu_1 \mid \nu_1 < \rho ) \sim 1/\rho$ as $\rho \to \infty$. Together with a similar argument for the second moment, we get
        \begin{equation}
    \label{eq:rho-moments-conditional}
    \Exp ( \nu_1 \mid \nu_1 < \rho ) \sim \frac{1}{\rho}, \text{ and } \Var ( \nu_1  \mid \nu_1 < \rho ) \sim \frac{1}{3\rho} .
\end{equation}
        By regeneration,
        $\sum_{n=1}^{\kappa_1} \nu_n$ is a sum of $\kappa_1$ independent random
        variables with the distribution of $\nu_1$ given $\nu_1 < \rho$,
        and  $\Pr(\kappa_1 =k)=\theta(\rho)(1-\theta(\rho))^k$ for $k\in \ZP$.
        Thus we 
        apply Lemma~\ref{lem:exponential} with $p_\rho \sim (4/\rho) \re^{-\pi^2 \rho/2}$ and $\mu_\rho \sim 1/\rho$ to conclude that
        $4 \re^{-\pi^2 \rho/2} \Lambda_1$ converges in distribution to $\cE$.
        The argument for $\Lambda_2$ is very similar: note that $\Pr ( \nu_1 > \rho - \rho^{-1/2} )$ also satisfies~\eqref{eq:rho-tail}.
\end{proof}

\begin{proof}[Proof of Theorem~\ref{thm:confined-space}]
We use a similar sandwiching argument to Lemma~\ref{lem:lambda-approx}.
Set 
\[ L_\rho := \sum_{n=1}^{\kappa_1+1} \ell_n \]
and define the event
\[ F_\rho := \{ \ell_n < \rho^{-1/2} \, \forall n \leq \kappa_1 + 1\}.\]
The first high-level excursion of duration exceeding $\rho$ forces
termination no later than $\rho$ units after its starting time, so
\begin{equation}
\label{eq:confined-bound-1}
\lambda_{1, \rho} \leq \sigma_{\kappa_1+1} + \rho = \Lambda_1 + L_\rho .\end{equation}
On $F_\rho$, for every $n \leq \kappa_2 \leq \kappa_1$,
\[ \gamma_n - \gamma_{n-1} = \ell_n + \nu_n < \rho^{-1/2} + (\rho - \rho^{-1/2} ) = \rho .\]
Starting from the boundary at time $\gamma_{n-1}$, the process thus survives until time $\gamma_n$. Inductively, it survives until $\gamma_{\kappa_2}$. Thus, on $F_\rho$,
\begin{equation}
\label{eq:confined-bound-2} \lambda_{1,\rho} \geq \gamma_{\kappa_2}+\rho \geq \rho + \sum_{n=1}^{\kappa_2} \nu_n = \Lambda_2 .\end{equation}
Let $T := \inf \{ t \geq 0 : | B_t | =1\}$ for standard Brownian motion $B_t$ with $B_0 =0$. By Brownian scaling, the conditional law of $\ell_n$ given $\cF_{\gamma_{n-1}}$ is the same as the law of $\rho^{-2} T$. Since $\Exp T =1$ and Lemma~\ref{lem: Lemma 2.1} gives $\Pr ( T \geq t ) \leq C \re^{-ct}$, we have
\[ \Exp \ell_n = \rho^{-2}, \text{ and } q_\rho := \Pr ( \ell_n \geq \rho^{-1/2} ) \leq C \re^{-c \rho^{3/2}} .\]
Consider the event
\[ A_n := \{ \kappa_1 + 1 \geq n \} = \{ \nu_j \leq \rho \text{ for all } 1 \leq j < n \} .\]
Then
\[ \Exp L_\rho = \sum_{n \geq 1} \Exp [ \ell_n \2{A_n} ] = \rho^{-2} \sum_{n \geq 1} (1 - \theta (\rho))^{n-1} = \frac{\rho^{-2}}{\theta(\rho)} .\]
Since $\theta(\rho) \sim (4/\rho) \re^{-\pi^2 \rho/2}$ we get
\[ \re^{-\pi^2 \rho/2} \Exp L_\rho \sim 1/(4\rho) \to 0 ,\]
so $\re^{-\pi^2 \rho/2} L_\rho \to 0$ in probability. Moreover,
\begin{align*}
\Pr ( F_\rho^\rc) & \leq \sum_{n \geq 1} \Pr ( A_n , \, \ell_n \geq \rho^{-1/2} ) \\
& = q_\rho \sum_{n \geq 1} ( 1 - \theta (\rho))^{n-1} \\
& = \frac{q_\rho}{\theta(\rho)} \to 0 ,\end{align*}
as $\rho \to \infty$. Consequently, for every $t \geq 0$,
\[ \Pr ( 4 \re^{-\pi^2 \rho/2} ( \Lambda_1 + L_\rho ) \leq t )
\leq \Pr ( 4 \re^{-\pi^2 \rho/2}  \lambda_{1,\rho} \leq t )
\leq  \Pr ( 4 \re^{-\pi^2 \rho/2}  \Lambda_{2} \leq t ) + \Pr ( F_\rho^\rc) .\]
Slutsky's theorem and the convergence results in Lemma~\ref{lem:lambda-asymptotics-confined} then yield $4 \re^{-\pi^2 \rho/2}  \lambda_{1,\rho}$ converges in distribution to $\cE$. By Brownian scaling, the claimed result follows.
\end{proof}

\subsection{The critical regime}
\label{sec:critical}

In this section we consider the case where $M, N \to \infty$ with $M/ N^2 \to \rho \in (0,\infty)$.
To state our main result, define the decreasing function $H$ by
\begin{equation}\label{eq:2.15}
    H(y):=\sum_{k=1}^\infty h_k(y), \quad h_k(y):=\exp\left\{-\frac{\pi^2(2k-1)^2y}{2}\right\}
\end{equation}
and for every $\rho>0$ and $s\in \R$ define,
\begin{equation}\label{eq:2.16}
    G(\rho,s):=\frac{s}{H(\rho)}\int_0^1e^{sv}\left(H(v\rho)-H(\rho)\right) \ud v,
\end{equation}
and for $\rho>0$ define $s_\rho:=\sup\{s>0:G(\rho,s)<1\}$. Then define, for $s<s_\rho$,
\begin{equation}\label{eq:2.17}
    \phi_\rho(s):=\frac{1}{1-G(\rho,s)}.
\end{equation}

  \begin{theorem}
    \label{thm:critical}
    Let $M, N \in (0,\infty)$ with $M/N^2 \to \rho \in (0,\infty)$. Let $X$ be reflected Brownian motion on $I_N$ with lifetime $\lambda$ as defined at~\eqref{eq:lambda-def}. 
        Then, as $M,N \to \infty$,                \begin{equation*}
            \frac{\lambda}{M} \tod 1+\xi_\rho, 
        \end{equation*}
        where  $\xi_\rho$ has moment generating function $\phi_\rho(s)$ from~\eqref{eq:2.17} for $s<s_\rho$.
    \end{theorem}
\begin{remarks}
    \label{rems:critical}
     \begin{myenumi}[label=(\alph*)]
\item
 If we use $\lambda_{N,M}$ to make the parameter dependence for the lifetime explicit, Brownian scaling and Theorem~\ref{thm:critical} show the exact distributional equality,  for every $\rho >0$,
 \begin{equation}
 \label{eq:critical-equality} \rho^{-1} \lambda_{1,\rho} \eqd 1+\xi_\rho .\end{equation}
 \item
    The distribution~$\xi_\rho$ first appeared, as far as we are aware,
    in the corresponding limit theorem for the critical case in the random walk model of~\cite{Wade_2023}: see Theorem~2.7 there.
    The continuous model studied here gives, via~\eqref{eq:critical-equality}, an exact probabilistic representation of the distribution that previously appeared only as a limit.
    \end{myenumi}
\end{remarks}

We use a similar high-level and low-level
excursion construction to Section~\ref{sec:BM-on-RP}. 
Recall that $X_0 = 0$. Define
$\gamma_0 := 0$ and then, recursively, for $n \in \N$,
\begin{align}
    \label{eq:sigma-def3}
\sigma_n & := \inf \{ t \geq \gamma_{n-1} : X_t  \in \{ 1 , N - 1 \}\} ; \\
    \label{eq:gamma-def3}
    \gamma_n & := \inf \{ t \geq \sigma_n : X_t  \in \{0,N\} \}.
\end{align}
Similarly to~\eqref{eq:kapp-def}--\eqref{eq:2.8}, define
\[    \kappa :=\max \{n\in \ZP: \gamma_n \leq \lambda\}, 
\]
and 
  \[ 
         \ell_n :=\sigma_n -\gamma_{n-1}  , \text{ and } 
        \nu_n :=\gamma_n -\sigma_n ,
\]
the durations of successive low-level and high-level excursions, respectively.
Similarly to~\eqref{eq:kappa-high}, define
\[ \kappa_1 := \min \{ n \in \N : \nu_n > M \} -1, \text{ and }
\theta (M) := \Pr ( \nu_1 > M).\]
Also, similarly to~\eqref{eq:kappa-high-new}, define, for $\eps>0$,
\[
\kappa_2 := \min \{ n \in \N : \nu_n > M - M^\eps \} - 1,
\]
and then, similarly to~\eqref{eq:lambda-high}, define
\begin{equation}
    \label{eq:lambda-high-critical}
    \Lambda_1 := M +\sum_{n=1}^{\kappa_1} \nu_n, \text{ and } \Lambda_2 := M + \sum_{n=1}^{\kappa_2} \nu_n.
\end{equation}
Then Lemma~\ref{lem:lambda-approx} applies essentially verbatim:

 \begin{lemma}
 Let $\eps \in (0,1)$. 
Define $\Lambda_1$ and $\Lambda_2$ as at~\eqref{eq:lambda-high-critical}, and consider event $E_M$ defined at~\eqref{eq:2.11} using the notation of this section. Then (i) $\lambda \leq \Lambda_1 + L_M$, and (ii) on $E_M$, $\lambda \geq \Lambda_2$.
    \end{lemma}


To find the asymptotics of $\Lambda_1$ and $\Lambda_2$ we will need the following result.

\begin{lemma}\label{lem:2.12}
    Let $H(y)$ be as defined at~\eqref{eq:2.15}. Consider a reflected Brownian motion on $I_N$. Then, for any $y_0 \in (0,\infty)$,
    \begin{equation*}
        \lim_{N\to\infty} \sup_{y \geq y_0} \left|\frac{N}{4}\Pr (\nu_1 >yN^2)-H(y)\right|=0.
    \end{equation*}
\end{lemma}
\begin{proof}
From~\eqref{eq:rho-full} we obtain  the formula
 \begin{equation*}
    \Pr ( \nu_1 >yN^2)=\frac{4}{\pi}\sum_{k=1}^\infty \frac{\sin\left(\frac{(2k-1)\pi}{N}\right)}{2k-1}\exp\left\{-\frac{(2k-1)^2\pi^2y}{2}\right\}.
\end{equation*}
 For convenience put $a_k = 2k-1$ and $c = \pi^2/2$. Let 
 \[ s_{N,k} := \frac{N}{\pi a_k} \sin \frac{\pi a_k}{N} .\]
 Then for every fixed $k \in \N$, $s_{N,k} \to 1$ as $N \to \infty$.
 Moreover, using $| \sin x | \leq x$ for $x \geq 0$, $| s_{N,k} | \leq 1$ so $| s_{N,k} -1 | \leq 2$.
 Hence, for fixed $y_0 >0$,
 \[ \sup_{y \geq y_0} \left| \frac{N}{4} \Pr ( \nu_1 > y N^2 ) - H(y) \right| \leq \sum_{k=1}^\infty | s_{N,k} - 1 | \re^{-c a_k^2 y_0} \to 0 ,\]
 as $N \to \infty$, by dominated convergence.
\end{proof}

The following is Lemma 4.8 from~\cite{Wade_2023}.

\begin{lemma}
\label{lem:4.8} Let $H(y)$ be as defined in \eqref{eq:2.15}. Then,
    \begin{enumerate}[label=(\roman*)]
        \item As $y\downarrow 0$, $H(y)\sqrt{8\pi y}\to 1$.
        \item As $y\to\infty$, $H(y) \sim e^{-\pi^2y/2}$
    \end{enumerate}
\end{lemma}

Now we can state the key results on asymptotics of $\nu_1$ that we need.

\begin{proposition}\label{prop:2.13}
    Suppose that 
     $M, N \to \infty$ with $M/ N^2 \to \rho \in (0,\infty)$.
Then,
    \begin{equation*}
        \theta(M) = \Pr ( \nu_1 > M) \sim \frac{4}{N} H(\rho).
    \end{equation*}
    Moreover, for any $s_0\in (0,\infty)$, as $M\to\infty$, uniformly for $s \in [-s_0,s_0]$ ,
    \begin{equation*}
        \Exp \left[ \re^{s\nu_1/N^2} \Bigmid \nu_1 \leq M \right]
       =
        1+\frac{4s}{N}\int_0^\rho \re^{sy}\left(H(y)-H(\rho)\right) \ud y + o (N^{-1}) .
    \end{equation*}
\end{proposition}
\begin{proof}
  Put
$r_N= M/ N^2$ and 
$Q_N(y)= (N/4) \Pr (\nu_1>yN^2)$.
Then Lemma~\ref{lem:2.12}
 and continuity of $H$ imply
$Q_N(r_N) \to H(\rho)$,
which proves the first assertion in the proposition.

For the second assertion, we also need control near zero. The series in the proof of Lemma~\ref{lem:2.12} gives
$0\le Q_N(y)\le H(y)$, for all $y>0$, 
because $|s_{N,k}|\le1$. Moreover, by
Lemma~\ref{lem:4.8}, $H(y)=O(y^{-1/2})$ as $y\downarrow0$, so $H$ is integrable on every bounded interval.
Then we obtain
\begin{align*}
\Exp \left[\re^{s\nu_1/N^2}\mid\nu_1\le M\right]-1
&=\frac{s}{1-\theta(M)}
\int_0^{r_N} \re^{sy}
\left[ \Pr (\nu_1>yN^2)-\theta(M)\right] \ud y\\
&=\frac{4s}{N(1-\theta(M))}
\int_0^{r_N} \re^{sy}\bigl[Q_N(y)-Q_N(r_N)\bigr] \ud y.
\end{align*}
Dominated convergence, using $Q_N(y)\le H(y)$, now yields
\[
\int_0^{r_N} \re^{sy}\bigl[Q_N(y)-Q_N(r_N)\bigr] \ud y \to 
\int_0^\rho \re^{sy}\bigl[H(y)-H(\rho)\bigr] \ud y,
\]
uniformly for $s\in[-s_0,s_0]$: on a common bounded integration interval, the exponential weights are uniformly bounded. Since $\theta(M)=O(N^{-1})$, the claimed result follows.
\end{proof}

\begin{lemma}\label{lem:2.15}
 Let $M, N \in (0,\infty)$ with $M / N^2 \to \rho \in (0,\infty)$. 
Define $\Lambda_1$ and $\Lambda_2$ as at~\eqref{eq:lambda-high-critical}.
  Then,  
  as $M \to \infty$, 
    \begin{equation*}
        \frac{\Lambda_1}{M} \tod 1+\xi_\rho, \text{ and } \frac{\Lambda_2}{M} \tod 1+\xi_\rho,
    \end{equation*}
    where  $\xi_\rho$ has moment generating function $\phi_\rho(s)$ from~\eqref{eq:2.17} for $s<s_\rho$.
\end{lemma}
\begin{proof}
Define
$S_M:=\sum_{n=1}^{\kappa_1}\nu_n$, $p_M:=\mathbb P(\nu_1>M)$, and
\[ \psi_M(s/M):=\mathbb E[\re^{s\nu_1/M}\mid\nu_1\le M].
\]
Then
\[
\mathbb E[\re^{sS_M/M}]
=\frac{p_M}{1-(1-p_M)\psi_M(s/M)}
\]
whenever $(1-p_M)\psi_M(s/M)<1$. Proposition~\ref{prop:2.13} gives
$p_M=(4/N) H(\rho)+o(N^{-1})$,
and
\[
\psi_M(s/M)
=1+\frac{4s}{N}
\int_0^1\re^{su}\bigl(H(\rho u)-H(\rho)\bigr)\ud u
+o(N^{-1}),
\]
using the fact that $M/N^2\to\rho$.
Consequently,
\[
1-(1-p_M)\psi_M(s/M)
=\frac4N H(\rho)\bigl(1-G(\rho,s)\bigr)+o(N^{-1}).
\]
For $0<s<s_\rho$, this denominator is eventually positive. Thus
\[
\mathbb E[ \re^{sS_M/M}]\longrightarrow
\frac1{1-G(\rho,s)}=\phi_\rho(s).
\]
The same calculation holds for negative $s$, so the moment-generating-function continuity theorem yields
\[
\frac{\Lambda_1}{M}=1+\frac{S_M}{M}
\tod 1+\xi_\rho.
\]
A similar argument applies for $\Lambda_2$.
\end{proof}

\begin{proof}[Proof of Theorem~\ref{thm:critical}.]
Essentially the same sandwiching argument as in the proof of Theorem~\ref{thm:half-line} applies to show that the convergence results for $\Lambda_1, \Lambda_2$ given by Lemma~\ref{lem:2.15} transfer to convergence of $\lambda$.
\end{proof}

\section{Multidimensional energy-constrained Brownian motion}
\label{sec:multidim}

\subsection{Overview}
\label{sec:multidim-overview}

In this section we make some initial steps into higher dimensions.
Section~\ref{sec:Energy constrained Brownian Motion The Confined Space regime} concerns the confined-space regime for the multidimensional case when the
region is a ball. Here we can use the radial description of multidimensional Brownian motion in terms of the Bessel process to carry out an analysis similar to the one-dimensional case. In Section~\ref{sec:universality} we make some comments and conjectures about general domains.

\subsection{The confined-space regime in a ball}
\label{sec:Energy constrained Brownian Motion The Confined Space regime}

Let $d \in \N$.
In this section we consider the domain $\cD_N := \{ x \in \R^d : \| x \| \leq N \}$,
the ball of radius $N$.
Let $\bX = (\bX_t)_{t \in \RP}$ denote reflected Brownian motion in $\cD_N$. 
For $t \in \RP$, define $e_t := t - \sup\{ s \in [0,t] : \bX_s \in \partial \cD_N\}$,
and again define the lifetime $\lambda$ by~\eqref{eq:lambda-def}. For simplicity we
throughout this section assume that
$\bX_0\in\partial D_N$,  equivalently $\|\bX_0\|=N$.

Our main result of this section is the asymptotic distribution of the lifetime in the confined-space regime, that is, when $M \gg N^2$. Let $J_\eta$ denote the Bessel function of the first kind of order $\eta$, and let $j_{\eta,1}, j_{\eta,2}, \ldots$ denote the increasing sequence of positive zeros of the Bessel function $J_\eta$.

\begin{theorem}
\label{thm:ball}
Let $d \in \N$.    As $M /N^2 \to \infty$,
    \[ \frac{d}{N^2} \exp \left\{ - j_{\frac{d}{2}-1,1}^2 \frac{M}{2N^2} \right\} \lambda \tod \cE ,\]
    where $\cE$ is a unit-mean exponential random variable.
\end{theorem}

\begin{remarks}
\label{rem:ball}
  \begin{myenumi}[label=(\alph*)]
\item   Taking $d=1$ and $N/2$ in place of $N$ (so we get reflected Brownian motion on the interval $[-N/2,N/2]$) we have $j_{-1/2,1} = \pi/2$ and the statement of Theorem~\ref{thm:ball}
becomes
\[ \frac{4}{N^2} \exp \left\{ - \pi^2 \frac{M}{2N^2} \right\} \lambda \tod \cE ,\]
recovering the statement of Theorem~\ref{thm:confined-space}. Another simple explicit value is when $d=3$, in which case $j_{1/2,1} = \pi$.
\item The   Bessel function of the first kind  $J_\eta (x)$ is given by~\cite[p.~360]{as} 
 \[
         J_\eta(x)=\sum_{k=0}^\infty \frac{(-1)^k}{k! \Gamma(\eta+k+1)}\left(\frac{x}{2}\right)^{2k+\eta},
   \]  and solves the Bessel equation $x^2\frac{\ud^2 y}{\ud x^2}+x\frac{\ud y}{\ud x}+(x^2-\eta^2)y=0$~\cite[p.~358]{as}. 
   \end{myenumi}
\end{remarks}

The norm process given by $R_t := \| \bX_t \|$ is a $d$-dimensional Bessel process on $[0,N]$ with reflection at $N$. To understand the lifetime of energy-constrained Brownian motion in $\cD_N$ it is enough to work with the one-dimensional process $R_t$,
since $e_t = t - \sup\{ s \in [0,t] : R_t = N\}$. We approach the analysis in a similar way to Section~\ref{sec:confined-space}, but using results on Bessel-process hitting times in place of hitting times for one-dimensional Brownian motion.

Making the dependence on $N,M$ explicit, Brownian scaling gives $\lambda_{N,M} \eqd N^2 \lambda_{1,\rho}$,
where $\rho = M/N^2$. Thus from this point onwards we will work in the unit ball with maximal energy capacity $\rho \to \infty$.
Let $T:= \inf \{ t \geq 0 : R_t =1\}$, and write $\Pr_x$ for the law of the Bessel process $R = (R_t)_{t \in \RP}$ started at $R_0 = x$. The next result concerns the tails of $T$.

\begin{lemma}
    \label{lem:Bessel-tail}
    Write $\eta = \frac{d}{2}-1$ and $j_k = j_{\eta,k}$ for short. 
    For $r \in (0,1)$ and $t >0$,
    \begin{equation}
    \label{eq:A-def} \Pr_r ( T >t) = \sum_{k=1}^\infty A_k (r) \re^{-j_k^2 t/2}, \text{ where } A_k (r) := \frac{2 r^{-\eta} J_\eta ( j_k r)}{j_k J_{\eta+1} (j_k)} .\end{equation}
    Moreover, there are constants depending only on $d$ such that, for $0 < \eps < 1/2$ and $t \geq 1$,
    \begin{align}
        \label{eq:T-bound-1}
\Pr_{1-\eps} ( T > t ) & \leq C \eps \re^{-j_1^2 t/2} ; \\
        \label{eq:T-bound-2}
      \Pr_{1-\eps} ( T > t ) & = ( 2 \eps + O (\eps^2) ) \re^{-j_1^2 t/2} + O ( \eps \re^{-j_2^2 t/2} ) .  
    \end{align}
\end{lemma}
\begin{proof}
    The explicit tail formula can be found in~\cite{Kent1980} (see also equation (2.8) of~\cite{hamana2012probabilitydistributionshittingtimes} with $a=r$ and $b=1$). We derive the uniform asymptotics as $r \uparrow 1$. It follows from Bessel recurrence relation
    \begin{equation}
    \label{eq:bessel-recurrence}
J_{\eta+1} (x) = \frac{\eta}{x} J_\eta (x) - J'_\eta (x) 
    \end{equation}
(see~\cite[p.~361]{as}) and    the simplicity of the zeros $j_k$~\cite[p.~370]{as} that $J_{\eta+1} (j_k) = -J'_\eta (j_k) \neq 0$. Also,  
\begin{equation}
    \label{eq:bessel-asympt}
    J_\eta (x) = \sqrt{\frac{2}{\pi x} } \left[ \cos \left( x - \frac{\eta \pi}{2} - \frac{\pi}{4} \right) + O (x^{-1} ) \right] 
\end{equation}
(see~\cite[p.~364]{as}) 
and (see~\cite[p.~371]{as})
\[
    j_k = \pi \left( k + \frac{\eta}{2} - \frac{1}{4} \right) + O (k^{-1}) .
\]
Consequently,
\begin{align*}
    J_{\eta+1} (j_k) = \sqrt{\frac{2}{\pi j_k
    } } \left[ \cos ( \pi (k-1) + O (k^{-1} ) ) + O( k^{-1} ) \right] ,
\end{align*}
which gives
\begin{equation}
    \label{eq:J-shifted}
    | J_{\eta+1} (j_k) | = \sqrt{\frac{2}{\pi j_k
    } } \left( 1 + O (k^{-1}) \right) .
\end{equation}
Also from~\eqref{eq:bessel-recurrence} we have
\begin{align*}
    \frac{\ud}{\ud r} \left[ r^{-\eta} J_\eta ( j_k r ) \right] & = - j_k r^{-\eta} J_{\eta+1} ( j_k r) ,
\end{align*}
and hence the derivative of $A_k(r)$ defined at~\eqref{eq:A-def} satisfies
\begin{equation}
    \label{eq:A-diff}
 A'_k(r) = - \frac{2 r^{-\eta} J_{\eta+1} (j_k r)}{J_{\eta+1} (j_k) } .
 \end{equation}
 By~\eqref{eq:J-shifted} there exist $k_0 \in \N$ and $c_0 >0$ such that $| J_{\eta+1} (j_k) | \geq c_0 j_k^{-1/2}$ for all $k \geq k_0$, while by~\eqref{eq:bessel-asympt} there exists $C< \infty$ such that, for all $r \in [1/2,1]$ and all $k$ sufficiently large, $| J_{\eta+1} (j_k r) | \leq C (j_k r)^{-1/2}$. Using these bounds in~\eqref{eq:A-diff} we see that $\sup_{1/2 \leq r \leq 1} | A_k' (r)|$ is bounded for $k$ sufficiently large, and hence (because $J_{\eta+1} (j_k) \neq 0$) for some $C< \infty$,
 \[ \sup_{k \in \N} \sup_{1/2 \leq r \leq 1 } | A_k'(r) | \leq C .\]
 Since $A_k (1) =0$, it follows that
 \begin{equation}
 \label{eq:A-near-1}
     | A_k (1-\eps) | \leq C \eps , \text{ for all } 0 < \eps \leq 1/2 \text{ and all } k \in \N.
 \end{equation}
For the first coefficient, \eqref{eq:A-diff} gives $A_1' (1) = -2$. The function $A_1$ has bounded second derivative on $[1/2,1]$ so Taylor's theorem gives $A_1 (1-\eps) = 2\eps + O (\eps^2)$ for $0 < \eps \leq 1/2$. For $m \in \{1,2\}$ define
\[ S_m := \sum_{k \geq m} \exp \left( - \frac{(j_k^2 - j_m^2)}{2} \right) ,\]
which has $S_m < \infty$ since $j_k \asymp k$. For $t \geq 1$,
\[ \sum_{k \geq m} \re^{-j_k^2 t/2} \leq S_m \re^{-j_m^2 t/2} .\]
Taking $m=1$ and using~\eqref{eq:A-near-1} we get
\[ \Pr_{1-\eps} (T > t ) \leq C \eps \sum_{k \geq 1} \re^{-j_k^2 t/2} \leq C S_1 \eps \re^{-j_1^2 t/2} ,\]
which gives~\eqref{eq:T-bound-1}. Similarly, using $m=2$,
\[ \left|  \Pr_{1-\eps} (T > t )  - A_1 (1-\eps) \re^{-j_1^2 t/2} \right| \leq C S_2 \eps \re^{-j_2^2 t/2} ,\]
and now use that $A_1 (1-\eps) = 2\eps + O (\eps^2)$ to get~\eqref{eq:T-bound-2}.
\end{proof}

Next we need results on the first two moments of $T$.

\begin{lemma}
\label{lem:T-moments}
    For $r \in [0,1]$, we have
    \[ \Exp_r T = \frac{1-r^2}{d}, \text{ and } \Exp_r ( T^2 ) = \frac{2(1-r^2)}{d^2} - \frac{1-r^4}{d(d+2)} .\]
    Consequently, as $\eps \downarrow 0$,
    \[ \Exp_{1-\eps} T \sim \frac{2\eps}{d}, \text{ and } \Exp_{1-\eps} ( T^2) \sim \frac{8\eps}{d^2 (d+2)} .\]
\end{lemma}
\begin{proof}
   The Laplace transform of $T$ is given by (see~\cite{gs} and equations (2.1)--(2.2) of~\cite{hamana2012probabilitydistributionshittingtimes})
   \[ \Phi_r (s) := \Exp_r \re^{-sT} = \frac{ r^{-\eta} I_\eta ( r \sqrt{2s} )}{I_\eta (\sqrt{2s})} ,\]
   where $I_\eta$ is the modified Bessel function of the first kind, defined by
   $I_\eta(x)=i^{-\eta}J_\eta(ix)$.
   
   Put $q = \eta+1= d/2$ and consider
   \begin{equation}
   \label{eq:F-series} F_r (s) = \sum_{k=0}^\infty \frac{ (r^2 s/2)^k}{k! (q)_k} = 1 + \frac{r^2}{d}s + \frac{r^4}{2d(d+2)} s^2 + O (s^3) ,\end{equation}
  as $s \to 0$, where $(q)_0 =1$ and $(q)_k = q (q+1) \cdots (q+k-1)$. Then $\Phi_r (s) = F_r(s)/ F_1 (s)$. Each $F_r$ is entire and $F_r (0)=1$, so $\log F_r$ is analytic near zero and, by~\eqref{eq:F-series},
   \[ \log F_r (s) = \frac{r^2}{d} s - \frac{r^4}{d^2 (d+2)} s^2 + O (s^3) .\]
   It follows that, as $s \to 0$,
   \begin{equation}
       \label{eq:log-Phi}
       \log \Phi_r (s) = - \frac{1-r^2}{d} s + \frac{1-r^4}{d^2 (d+2)} s^2 + O (s^3) .
   \end{equation}
   Then $\Exp_r T = - \Phi'_r (0) = - (\log \Phi_r)' (0)$ gives the claimed formula
   for $\Exp_r T$, and $\Var_r T = \Phi''_r (0) - ( \Phi_r'(0))^2 = (\log \Phi_r)'' (0)$
   gives
   \[ \Var_r T = \frac{2 (1-r^4)}{d^2 (d+2)} .\]
   The formula for $\Exp_r (T^2)$ follows. The asymptotics for $r = 1-\eps$, $\eps \downarrow 0$, follow
   from the expansions $1-r^2 = 2\eps + O (\eps^2)$ and $1-r^4 = 4\eps + O (\eps^2)$.
\end{proof}

We now introduce notation for excursions, in a similar way to Section~\ref{sec:confined-space}. Take $\rho >2$ and $\eps = \rho^{-1}$. Starting with $R_0 = 1$ define $\gamma_0 := 0$ and then, recursively,
for $n \in \N$,
\[ \sigma_n := \inf \{ t \geq \gamma_{n-1} : R_t = 1-\eps\}, \text{ and } \gamma_n := \inf \{ t \geq \sigma_n : R_t =1\} .\]
For $n \in \N$ let
\[ \ell_n := \sigma_n -\gamma_{n-1}, \text{ and } \nu_n := \gamma_n - \sigma_n ,\]
the durations of low-level and high-level excursions, respectively.

Set $\delta = \rho^{-1/2}$, $h_1 = \rho$ and $h_2 = \rho -\delta$, and define, for $i \in \{1,2\}$,
\[ \kappa_i := \min \{ n\in \N : \nu_n > h_i \} -1, \text{ and } \Lambda_i := \rho + \sum_{n=1}^{\kappa_i} \nu_n .\]
The following result gives the key asymptotic estimates for high-level excursions.

\begin{lemma}
\label{lem:ball-asymptotics}
    For $i \in \{1,2\}$, as $\rho \to \infty$,
    \begin{align*}
        \Pr ( \nu_1 > h_i ) & \sim \frac{2}{\rho} \re^{-j_1^2 \rho/2} ; \\
        \Exp ( \nu_1 \mid \nu_1 \leq h_i) & \sim \frac{2}{d\rho} ; \\
        \Exp ( \nu_1^2 \mid \nu_1 \leq h_i) & \sim \frac{8}{d^2(d+2)\rho} .
    \end{align*}
\end{lemma}
\begin{proof}
    First note that $\Pr ( \nu_1 > h_i) = \Pr_{1-\eps} ( T > h_i)$, and then using~\eqref{eq:T-bound-2} we get the claimed asymptotics. For $\beta \in \{1,2\}$ and $h \geq 1$, we have from~\eqref{eq:T-bound-1} that
    \begin{align*}
        \Exp [ \nu_1^\beta \1 { \nu_1 \geq h} ] & = h^\beta \Pr ( \nu_1 > h) + \beta \int_h^\infty t^{\beta-1} \Pr ( \nu_1 > t ) \ud t \\
        & \leq \frac{C}{\rho} \left( h^\beta \re^{-j_1^2 h/2} + \int_h^\infty t^{\beta-1} \re^{-j_1^2 t/2} \ud t \right).
    \end{align*}
    The integral here can be evaluated to give
    \begin{align*}
         \int_h^\infty  \re^{-j_1^2 t/2} \ud t & = \frac{2}{j_1^2} \re^{-j_1^2 h/2} ; \\
         \int_h^\infty t \re^{-j_1^2 t/2} \ud t & = \left( \frac{2h}{j_1^2} + \frac{4}{j_1^4} \right) 
         \re^{-j_1^2 h/2} .
    \end{align*}
    It follows that $\Exp [ \nu_1^\beta \1{ \nu_1 \geq h_i } ] \leq C \rho^{\beta-1}  \re^{-j_1^2 h_i/2} = o(\rho^{-1})$. Then writing
    \[ \Exp ( \nu_1^\beta \mid \nu_1 \leq h_i ) = \frac{ \Exp ( \nu_1^\beta) - \Exp [ \nu_1^\beta \1{ \nu_1 > h_i } ]}{\Pr ( \nu_1 \leq h_i )} ,\]
    using that $\Exp (\nu_1^\beta) = \Exp_{1-\eps} ( T^\beta)$ and applying the
    asymptotics from Lemma~\ref{lem:T-moments} gives the claimed results.
\end{proof}

\begin{lemma}
\label{lem:Lambda-ball}
    For each $i \in \{1,2\}$,
    \[ d \re^{-j_1^2 \rho /2} \Lambda_i \tod \cE ,\]
    where $\cE$ is a unit-mean exponential random variable.
\end{lemma}
\begin{proof}
    Conditional on $\kappa_i =k$, the first $k$ high-level excursions are independent with the law of $\nu_1$ conditional on $\nu_1 \leq h_i$. Thus $\sum_{n=1}^{\kappa_i} \nu_n$ has the same distribution as $Z_\rho = \sum_{n=1}^{K_\rho} Y_{\rho,n}$ where the $K_\rho$ and $Y_{\rho,n}$ are independent, $\Pr ( K_\rho = k) = p_\rho (1-p_\rho)^k$ with $p_\rho = \Pr ( \nu_1 > h_i)$, and the $Y_{\rho,n}$ are i.i.d.~with the law of  $\nu_1$ conditional on $\nu_1 \leq h_i$. Applying Lemma~\ref{lem:exponential}
    with the estimates from Lemma~\ref{lem:ball-asymptotics} yields the result.
\end{proof}

The next result estimates the durations of low-level excursions, playing a similar role to Lemma~\ref{lem: Lemma 2.1}.

\begin{lemma}
    \label{lem:bessel-low}
    There exist $C<\infty$ and $c>0$ such that, for all $\eps \in (0,1/2]$ and all $t \geq 0$, \[ \Pr (\ell_1 \geq t ) \leq C \re^{-ct/\eps^2}. \]
\end{lemma}
\begin{proof}
    Suppose that $R_0 = y \in (1-\eps,1]$ and recall that $\sigma_1 = \inf \{ t \geq 0 : R_t = 1-\eps\}$. Up until time $\sigma_1$, $R_t$ satisfies the SDE
    \[ \ud R_t = \frac{d-1}{2R_t} \ud t + \ud W_t - \ud L_t ,\]
    where the local-time term $L_t$, governing reflection at~$1$, increases only when $R_t = 1$. Before time $\sigma_1$, $R_t > 1-\eps \geq 1/2$, so
    it follows that 
    \[ R_t \leq y + (d-1) t + W_t, \text{ on } \{ t \leq \sigma_1 \} .\]
   On the event $W_h < 1 - \eps - (d-1)h -y$, if $h \leq \sigma_1$ the above relation implies that $R_h < 1-\eps$, a contradiction. Thus
   \[ \{ W_h <   - \eps  - (d-1)h \} \subseteq \{ W_h < 1 - \eps - (d-1)h -y \} \subseteq \{ \sigma_1 \leq h \} .\]
   In particular, with $h = \eps^2$ we get
   \[ \Pr ( W_h <   - \eps  - (d-1) \eps^2 )
   = \Pr ( W_1 < -1 - (d-1)\eps ) \geq \Pr ( W_1 < -1 - (d-1)/2 ) =: p_d ,\]
   where $p_d >0$. 
   By the Markov property at time $t$ it follows that
   \[ \Pr ( \sigma_1 \leq t + \eps^2 \mid \cF_t ) \geq p_d , \text{ on } \{ \sigma_1 \geq t \} ,\]
   and the result follows by iterating this bound.
\end{proof}

Now we can complete the proof of Theorem~\ref{thm:ball}.

\begin{proof}[Proof of Theorem~\ref{thm:ball}]
Recall that $\eps = \rho^{-1}$ and $\delta = \rho^{-1/2}$. Then, by Lemma~\ref{lem:bessel-low}, $\Pr ( \ell_1 \geq \delta) \leq C \re^{-c \rho^{3/2}}$ and $\Exp \ell_1 \leq C \eps^2 = O (\rho^{-2})$. Put
\[ L_\rho := \sum_{n=1}^{\kappa_1+1} \ell_n, \text{ and } F_\rho := 
\{ \ell_n < \delta \text{ for all } 1 \leq n \leq \kappa_1 +1 \}.\]
Let $A_n := \{ \kappa_1+1 \geq n\} = \{ \nu_j \leq \rho \text{ for all } j < n\}$. Then
\[ \Exp L_\rho = \sum_{n\geq 1} \Exp ( \ell_n \2{A_n} ) = \frac{\Exp \ell_1}{\Pr (\nu_1 \geq h_1)} ,\]
which together with Lemma~\ref{lem:ball-asymptotics} shows that
\begin{equation}
\label{eq:L-rho-small}
\re^{-j_1^2 \rho/2} \Exp L_\rho = O (\rho^{-1} ) \to 0 .\end{equation}
Similarly,
\begin{equation}
\label{eq:F-rho-small} 
\Pr ( F_\rho^\rc ) \leq \sum_{n \geq 1} \Pr ( A_n \cap \{ \ell_n \geq \delta\} ) = \frac{\Pr (\ell_1 \geq \delta )}{\Pr ( \nu_1 \geq h_1)} \to 0.\end{equation}
Then, similarly to previous arguments,
$\lambda \leq \sigma_{\kappa_1+1} + \rho = \Lambda_1 + L_\rho$.
On $F_\rho$, for $n \leq \kappa_2 \leq \kappa_1$, $\gamma_n - \gamma_{n-1} = \ell_n + \nu_n < \delta + \rho - \delta = \rho$, so the process survives until at least $\gamma_{\kappa_2}+\rho$, i.e.,
\[ \lambda \geq \gamma_{\kappa_2} + \rho = \rho + \sum_{n=1}^{\kappa_2} ( \ell_n + \nu_n) \geq \Lambda_2 \text{ on } F_\rho .\]
Write $c_\rho = d \re^{-j_1^2 \rho/2}$. Then we have
\[ \Pr ( c_\rho ( \Lambda_1 + L_\rho ) \leq t ) \leq \Pr ( c_\rho \lambda \leq t)
\leq \Pr ( c_\rho \Lambda_2 \leq t) + \Pr ( F_\rho^\rc ) .\]
This sandwiching argument together with Lemma~\ref{lem:Lambda-ball} and the estimates in~\eqref{eq:L-rho-small} and~\eqref{eq:F-rho-small} shows that $c_\rho \lambda$ converges in distribution to $\cE$, as claimed.
\end{proof}

\subsection{Some hints of universality}
\label{sec:universality}

On the basis of the results so far, we state in this section some conjectures for the case of general domains in general dimensions. Take $\cD$ to be a fixed smooth, connected, bounded domain in $\R^d$, and scale by parameter $N$ to get scaled domain $N \cD$. Consider reflected Brownian motion in $N \cD$ and define the lifetime $\lambda$ to be the first time that the elapsed time since the most recent visit to the boundary reaches parameter $M$. Assume that the process starts on the boundary.

In the meagre-capacity regime $M/N^2 \to 0$ we expect the Darling--Mandelbrot limit from Theorem~\ref{thm:meagre} is universal:
\begin{conjecture}
    \label{conj:meagre}
    Let $M, N \in (0,\infty)$ with $M/ N^2  \to 0$. For the case of general multidimensional domains described above, 
        for $\xi\sim \DM(1/2)$, 
        \[ \frac{\lambda}{M} \tod 1 +\xi.\]
\end{conjecture}
The intuition is that in the meagre-capacity regime it is unlikely for the process to visit distant sections of the boundary, and locally the process should be well-approximated by reflected Brownian motion on the half-space, which reduces to the case of $\RP$ analysed in Section~\ref{sec:BM-on-RP}.

We turn to the confined-space regime.
Let $\mu_1:=\mu_1(\cD)$ be the principal Dirichlet eigenvalue of $-\frac12\Delta$ on $\cD$, and let $\phi_1>0$ be a corresponding eigenfunction:
\[
-\frac12\Delta\phi_1=\mu_1\phi_1, \text{ in }\cD,
\qquad
\phi_1=0, \text{ on }\partial \cD.
\]
Define the  constant
\[
 {
\kappa_\cD
:=
\mu_1
\frac{\left(\displaystyle\int_\cD\phi_1(x) \ud x\right)^2}
{|\cD|\displaystyle\int_\cD\phi_1(x)^2 \ud x}.
}
\]

Here is the conjecture in the confined-space regime. 
\begin{conjecture}
    \label{conj:confined}
    Let $M, N \in (0,\infty)$ with $M/ N^2  \to \infty$. For the case of general multidimensional domains described above, 
        for $\cE$ a unit-mean exponential,  
\begin{equation}
    \label{eq:conjecture}
 {
\frac{\kappa_\cD}{N^2}
\exp\left\{-\mu_1 \frac{M}{N^2}\right\}
\lambda
\tod 
\cE.
}
\end{equation}
\end{conjecture}

In the critical regime, the intuition from Theorem~\ref{thm:critical} is that 
spectral information beyond the principal eigenvalue appears in the limit.
We do not state a specific conjecture for the general case.

To finish this section, we give some justification for Conjecture~\ref{conj:confined}. By Brownian scaling, it suffices to work in the fixed domain $\cD$, with energy capacity $\rho=M/N^2\to\infty$. Termination occurs on the first excursion from the boundary whose duration reaches $\rho$. 
In what follows, $\bX$ denotes normally reflected Brownian motion in $\cD$.
We use the local-time normalization employed by Hsu \cite{hsu-1986}:
\begin{equation}\label{eq:hsusde}
 \ud \bX_t=\ud \bB_t+\tfrac12 n(\bX_t)\ud  L_t,
\end{equation}
where $\bB$ is $d$-dimensional Brownian motion, $n$ is the  {inward} unit normal, and local-time term $L$
charges only boundary-contact times. Keeping track of the factor $1/2$ in~\eqref{eq:hsusde} is important for the constant in
\eqref{eq:conjecture}.  
Under this convention, \(L\) has Revuz measure \(\sigma\) (the surface area) relative to Lebesgue measure on \(D\); equivalently, its stationary mean occupation identity is
\begin{equation}\label{eq:localtime-intensity}
 \Exp_\pi \left[\int_0^T f(\bX_s)\ud  L_s\right]
 =\frac{T}{|\cD|}\int_{\partial\cD}f(z) \ud \sigma(z)
\end{equation}
for every bounded Borel measurable $f:\partial D\to\mathbb R$
(see~\cite{hsu-1986} equations (2.3)--(2.4), and 
recall that the invariant probability measure of reflected
Brownian motion is the normalized Lebesgue measure $\pi=|\cD|^{-1}\ud x$).

Let $\tau_D=\inf\{t\geq0:\bB_t\notin \cD\}$ for Brownian motion killed at
$\partial\cD$, and put
\begin{equation}\label{eq:u}
 u(t,x):=\Pr_x(\tau_\cD>t),\qquad t>0,\quad x\in \cD.
\end{equation}
Write $Q^z$ for the excursion measure from $z\in\partial\cD$ in Hsu's
normalization, and $\zeta$ for the lifetime of a single excursion. The
absolute excursion law in Hsu's Theorem~5.1~\cite{hsu-1986}, together with the normal
derivative formula for the killed heat kernel, gives
\begin{equation}\label{eq:exc-tail}
 Q^z(\zeta>t)
 =\frac12\,\partial_n u(t,z),\qquad t>0,
\end{equation}
where $\partial_n$ is the inward normal derivative. Indeed, if
$p_\cD(t,x,y)$ is the killed heat kernel, the excursion entrance density in
\cite[\S\S2,5]{hsu-1986} is
$\frac12\partial_{n_z}p_\cD(t,z,y)$, and integration over $y\in \cD$
yields \eqref{eq:exc-tail}.

Hsu's \cite[Theorem~5.1]{hsu-1986} compensator identifies the rate of
excursions conditional on the boundary trace process
$\xi_s=X_{L^{-1}(s)}$. Thus, on the local-time scale, the
predictable intensity of excursions of duration exceeding $t$ is
$Q^{\xi_s}(\zeta>t)\ud s$. In dimension $d\geq2$ this is generally
\emph{not} an ordinary Poisson point process with a deterministic rate:
the excursion law depends on the current boundary point. For a modern
boundary-process and exit-system framework see
\cite[Chapter~5, especially \S5.7]{chen-fukushima}.
The construction and basic potential theory of reflecting Brownian
motion and boundary local time under weaker geometric assumptions are
studied by Bass and Hsu \cite{bass-hsu}. 

Combining the excursion compensator with the stationary boundary
local-time intensity \eqref{eq:localtime-intensity} suggests the following
\emph{stationary mean rate} of excursions with lifetime greater than $t$:
\begin{equation}\label{eq:rate}
 r_\cD(t):=\frac1{|\cD|}\int_{\partial\cD}Q^z(\zeta>t)\ud\sigma(z)
 =\frac1{2|\cD|}\int_{\partial\cD}\partial_n u(t,z)\ud\sigma(z).
\end{equation}
This rate counts the starts of long excursions per unit of ordinary time
in stationarity.

Define the Dirichlet heat content 
by
\begin{equation}\label{eq:heatcontent}
 H_\cD(t):=\int_\cD u(t,x)\ud x
 =\int_\cD\Pr_x(\tau_\cD>t)\ud x.
\end{equation}
The function $u$ solves
\[
 \partial_tu=\tfrac12\Delta u\quad\text{in }\cD,
 \qquad u(t,z)=0\quad(z\in\partial\cD),
 \qquad u(0,x)=1\quad(x\in \cD).
\]
For $t>0$, integration of the heat equation and the divergence theorem
give (recalling that $n$ points  {inwards})
\begin{equation}\label{eq:heat-identity}
 H_\cD'(t)=\frac12\int_\cD\Delta u(t,x)\ud x
 =-\frac12\int_{\partial\cD}\partial_nu(t,z)\ud\sigma(z).
\end{equation}
Consequently the mean long-excursion rate has the particularly simple form
\begin{equation}\label{eq:rateheat}
 {\qquad r_\cD(t)=-\frac{H_\cD'(t)}{|\cD|}.\qquad}
\end{equation}

To extract its large-$t$ asymptotics, let
$\{(\mu_k,\phi_k):k\geq1\}$ be an orthonormal Dirichlet
eigenbasis for $-\tfrac12\Delta$ in $L^2(\cD)$, with
$0<\mu_1<\mu_2\leq\cdots$. The standard killed-semigroup
expansion (see, for example, Davies \cite{davies}) yields
\begin{equation}\label{eq:heatexpansion}
 H_\cD(t)
 =\langle 1,\re^{t\Delta_\cD/2}1\rangle_{L^2(\cD)}
 =\sum_{k=1}^\infty \re^{-\mu_k t}
   \left(\int_\cD\phi_k(x)\ud x\right)^2.
\end{equation}
Termwise differentiation is legitimate for fixed $t>0$; hence
\begin{equation}\label{eq:rateexpansion}
 r_\cD(t)=\frac1{|\cD|}\sum_{k=1}^\infty
 \mu_k \re^{-\mu_k t}
 \left(\int_\cD\phi_k(x)\ud x\right)^2.
\end{equation}
Since $\cD$ is connected, $\mu_1$ is simple and $\phi_1$ can be
chosen strictly positive. The first term in \eqref{eq:rateexpansion}
therefore dominates, giving
\begin{equation}\label{eq:asymptoticrate}
 r_\cD(t)\sim\kappa_\cD \re^{-\mu_1t}
 \qquad(t\to\infty), 
\end{equation}
with $\kappa_D$ as defined above. This establishes the proposed
asymptotic \emph{mean occurrence rate} of long excursions, subject to the
normalizations above.

The reflected process in $D$ is ergodic and has a finite relaxation
time, whereas \eqref{eq:asymptoticrate} predicts a typical waiting time
of order $\re^{\mu_1\rho}$ before an excursion longer than $\rho$
starts. The two scales are widely separated as $\rho\to\infty$.
Consequently one expects the boundary trace, and the ordinary-time
process, to have mixed between successive opportunities for a rare
long excursion.
This suggests that, after rescaling ordinary time by $r_D(\rho)$,
the point process of \emph{starts} of excursions whose duration reaches
$\rho$ should converge to a rate-one Poisson process. Its first point
would then yield 
an exponentially distributed waiting time. Once such an excursion
starts, the energy-constrained lifetime occurs exactly $\rho$ units
later. Because
$ \rho\,r_\cD(\rho)\to 0$,
this final additive delay is negligible on the proposed scale. Together
with \eqref{eq:asymptoticrate}, the prediction is
\begin{equation}\label{eq:scaled}
 \kappa_\cD \re^{-\mu_1\rho} \lambda
 \tod \cE,
\end{equation}
which is equivalent to \eqref{eq:conjecture}.

Finally, we check the conjecture against the interval and ball results that we proved above.

\paragraph{The unit interval.}
For $\cD=(0,1)$ take
$\mu_1=\pi^2/2$ and $\phi_1(x)=\sin(\pi x)$. Then
\[
 \int_0^1\phi_1(x)\ud x=\frac2\pi,
 \qquad \int_0^1\phi_1(x)^2\ud x=\frac12,
\]
so that
\begin{equation}\label{eq:interval}
 \kappa_{\cD}
 =\frac{\pi^2}{2}\frac{(2/\pi)^2}{1/2}=4.
\end{equation}
Thus \eqref{eq:conjecture} reads
\[
 \frac4{N^2}\exp\!\left\{-\frac{\pi^2M}{2N^2}\right\}
 \lambda \tod \cE,
\]
precisely the one-dimensional confined-space normalization in
Theorem~\ref{thm:confined-space}.

\paragraph{The unit ball.}
Let $\cD=B(0,1)\subset\R^d$, set
$\nu=d/2-1$ and $j=j_{\nu,1}$, the first positive zero of the Bessel
function $J_\nu$. Then separation of variables (see e.g.~\cite[Ch.\ V]{courant}) gives the radial equation
\[
\phi''+\frac{d-1}{r}\phi'+2\mu\phi=0;
\]
Putting \(\phi(r)=r^{-(d/2-1)}v(\sqrt{2\mu}\,r)\) gives Bessel’s equation of order \(d/2-1\).
Thus we have 
\[
 \mu_1=\frac{j^2}{2},\qquad
 \phi_1(x)=\lVert x\rVert^{-\nu}J_\nu(j\lVert x\rVert).
\]
Writing $S=|\partial B(0,1)|=d|B(0,1)|$, standard Bessel identities 
(see~equations
(11.3.20), (11.4.5), and (9.1.27) of~\cite{as}) 
yield
\begin{align*}
 \int_\cD\phi_1(x)\ud x
  &=S\int_0^1r^{\nu+1}J_\nu(jr)\ud r
    =\frac{S}{j}J_{\nu+1}(j),\\
 \int_\cD\phi_1(x)^2\ud x
  &=S\int_0^1rJ_\nu(jr)^2\ud r
    =\frac S2J_{\nu+1}(j)^2.
\end{align*}
It follows that
\begin{equation}\label{eq:ball}
 \kappa_\cD
 =\frac{j^2}{2|\cD|}\,
   \frac{S^2J_{\nu+1}(j)^2/j^2}
        {(S/2)J_{\nu+1}(j)^2}
 =\frac{S}{|\cD|}=d.
\end{equation}
Thus the conjecture gives  
\[
 \frac{d}{N^2}
 \exp\!\left\{-\frac{j_{d/2-1,1}^2 M}{2N^2}\right\}
 \lambda \tod \cE,
\]
matching the ball result in Theorem~\ref{thm:ball}.

\section*{Acknowledgements}
\addcontentsline{toc}{section}{Acknowledgements}

The research for this note was undertaken during an undergraduate summer project at Durham University.  ChatGPT helped with checking earlier drafts of the paper, suggesting improvements,  and, especially,  with formulating Conjecture~\ref{conj:confined} in Section~\ref{sec:universality}.

\printbibliography

\end{document}